\documentclass[11pt]{article}

\RequirePackage{amsthm,amsmath,amsmath,amsfonts,amssymb}
\RequirePackage[colorlinks,citecolor=blue,urlcolor=blue]{hyperref}

\RequirePackage{tablefootnote}

\RequirePackage{graphicx}

\usepackage{fullpage}

\renewcommand{\bar}{\overline} 

\newcommand{\bE}{{\mathbb{E}}}

\newcommand{\eqdef}{\stackrel{\textrm{def}}{=}}

\usepackage{amssymb}

\usepackage{enumitem} 
\setitemize{itemsep=0mm, leftmargin=5mm, topsep=0mm}
\setenumerate{itemsep=0mm, leftmargin=5mm, topsep=0mm}

\newtheorem{theorem}{Theorem}
\newtheorem{lemma}{Lemma}
\newtheorem{corollary}{Corollary}

\newtheorem{proposition}{Proposition}

\theoremstyle{definition}

\allowdisplaybreaks[4]

\begin{document}




\title{On Asymptotically Tight Tail Bounds for Sums of Geometric and Exponential Random Variables}


\author{
Yaonan Jin \thanks{Department of IEDA, Hong Kong University of Science and Technology. Email: \tt{yjinan@connect.ust.hk}.} \and
Yingkai Li \thanks{Department of Computer Science, Northwestern University. Email: \tt{yingkai.li@u.northwestern.edu}.} \and
Yining Wang \thanks{Machine Learning Department, Carnegie Mellon University. Email: \tt{yiningwa@cs.cmu.edu}.} \and
Yuan Zhou \thanks{Computer Science Department, Indiana University at Bloomington; Department of Industrial and Enterprise Systems Engineering, University of Illinois at Urbana-Champaign. Email: \tt{yuanz@illinois.edu}.}
}










\maketitle
\begin{abstract}
In this note we prove bounds on the upper and lower probability tails of sums of independent geometric or exponentially distributed
random variables.
We also prove negative results showing that our established tail bounds are asymptotically tight.
\end{abstract}

\section{Introduction}

Consider independent and identically distributed random variables $X_1,\cdots,X_n$
from either the \emph{geometric distribution} (the failure model)
\begin{equation}
\Pr[X=k] = p(1-p)^k, \;\;\;\;\;\;k=0,1,2,\cdots
\label{eq:geometric}
\end{equation}
or the \emph{exponential distribution}
\begin{equation}
p(y)=\rho e^{-\rho y}, \;\;\;\;\;\; y\in[0,\infty).
\label{eq:exponential}
\end{equation}

Let $\bar X_n:=(X_1+\cdots+X_n)/n$  be the normalized sum of i.i.d.~geometric or exponential random variables,
and $\mu := \mathbb EX_1 = (1-p)/p$ for Eq.~(\ref{eq:geometric}) or $\mu :=1/\rho$ for Eq.~(\ref{eq:exponential}) be the mean of each random variable.
The objective of this paper is to characterize the \emph{tail probabilities} $\Pr[\bar X_n\geq \lambda \mu]$ for $\lambda\in(1,\infty)$ or
$\Pr[\bar X_n\leq \lambda \mu]$ for $\lambda\in(0,1)$.
Such tail bounds are important in statistics, computer science and operations research, 
and have recently found interesting applications in assortment selection problems in operations management \cite{agrawal2017mnl,agrawal2017thompson,chen2018dynamic}.

Our main result can be summarized in the following theorem:
\begin{theorem}
For any $\lambda,\mu>0$, define $H(\lambda,\mu):=\mu\lambda\ln\lambda - (1+\mu\lambda)\ln\left(\frac{1+\mu\lambda}{1+\mu}\right)$ and
$G(\lambda):=\lambda-1-\ln\lambda$. 
Then for geometric random variables $X_1,\cdots,X_n$, 
\begin{align*}
\Pr[\bar X_n \geq \lambda\mu] &\leq \exp\left\{-n\cdot H(\lambda,\mu)\right\} \;\;\;\;\;\lambda\in(1,\infty); & \text{(upper tail)}\\
\Pr[\bar X_n \leq  \lambda\mu] &\leq \exp\left\{-n\cdot H(\lambda,\mu)\right\} \;\;\;\;\;\lambda\in(0,1).& \text{(lower tail)}\\
\end{align*} 
In addition, for exponential random variables $Y_1,\cdots,Y_n$, we have 
\begin{align*}
\Pr[\bar Y_n \geq \lambda\mu] &\leq \exp\left\{-n\cdot G(\lambda)\right\} \;\;\;\;\;\lambda\in(1,\infty); & \text{(upper tail)}\\
\Pr[\bar Y_n \leq  \lambda\mu] &\leq \exp\left\{-n\cdot G(\lambda)\right\} \;\;\;\;\;\lambda\in(0,1).& \text{(lower tail)}\\
\end{align*} 
\label{thm:main-positive}
\end{theorem}
\vspace{-0.5in}

 Note that for all $\mu,\lambda>0$, $H(\lambda,\mu)\geq 0$ and $G(\lambda)\geq 0$, with equality holds if and only if $\lambda=1$.
 This means both upper and lower tails of $\bar X_n$ or $\bar Y_n$ decay exponentially fast as $\exp\{-\Omega( n)\}$ provided that $\lambda\neq 1$.
 
Theorem \ref{thm:main-positive} is proved by careful applications of the \emph{Chernoff bound}, which is given in Sec.~\ref{sec:proof-positive}.
We also give several approximations of $H(\lambda,\mu)$ in Sec.~\ref{sec:approx}, 
making it easier for practical usage and also to compare against existing tail bounds \cite{janson2018tail,agrawal2017mnl}.

To understand the tightness of Theorem \ref{thm:main-positive}, and especially $H(\lambda,\mu)$ and $G(\lambda)$,
 we state the following result: 
 
\begin{corollary}
Let $\lambda,\mu>0$, $\lambda\neq 1$ be fixed and not changing with $n$.
For geometric random variables $X_1,\cdots,X_n,\cdots$, we have
\begin{align*}
&\lim_{n\to\infty}-{\ln \Pr[\bar X_n\geq \lambda\mu]}/[{n H(\lambda,\mu)} ]= 1\;\;\;\;\;\;\lambda\in(1,\infty);& \text{(upper tail)}\\
&\lim_{n\to\infty} -{\ln \Pr[\bar X_n\leq \lambda\mu]}/[{n H(\lambda,\mu)}]= 1\;\;\;\;\;\;\lambda\in(0,1).& \text{(lower tail)}
\end{align*}
Similarly, for exponential random variables $Y_1,\cdots,Y_n,\cdots$, we have 
\begin{align*}
&\lim_{n\to\infty} -{\ln \Pr[\bar Y_n\geq \lambda\mu] }/[{nG(\lambda)}]= 1\;\;\;\;\;\;\lambda\in(1,\infty);& \text{(upper tail)}\\
&\lim_{n\to\infty}  -{\ln \Pr[\bar Y_n\leq \lambda\mu] }/[{nG(\lambda)}] = 1\;\;\;\;\;\;\lambda\in(0,1).& \text{(lower tail)}
\end{align*}
\label{cor:negative}
\end{corollary}

Corollary \ref{cor:negative} is a simplified statement of Theorems \ref{thm:negative-geometric} and \ref{thm:negative-exponential},
both of which are stated and proved in Sec.~\ref{sec:proof-negative}.
It shows that the leading terms in the exponents of the tail bounds in Theorem \ref{thm:main-positive},
namely $-n\cdot H(\lambda,\mu)$ and $-n\cdot G(\lambda)$, are asymptotically tight as $n\to\infty$.
Actually, in Theorems \ref{thm:negative-geometric} and \ref{thm:negative-exponential} we show that the remainder terms are on the order of $O(\log n)$,
which is considerably smaller than the $-n\cdot H(\lambda,\mu)$ and $-n\cdot G(\lambda)$ leading terms.

Corollary \ref{cor:negative} is proved by binomial or Poisson counting process characterizations of $\bar X_n$ and $\bar Y_n$,
for which \emph{exact} tail probabilities are known via combinations of binomial coefficients.
Afterwards, Stirling's approximation is applied to derive asymptotic expressions of the tail probabilities.
The complete proof is given in Sec.~\ref{sec:proof-negative}.

\section{Approximations of $H(\lambda,\mu)$, and comparisons}\label{sec:approx}

In this section we give several approximations of $H(\lambda,\mu)$ with simpler forms.
We also compare our result with existing tail bounds for sums of geometric random variables,
mostly from \cite{janson2018tail} and \cite{agrawal2017mnl}, showing our bound is tighter in several cases and easier to use overall.

\begin{proposition}
The function $H(\lambda,\mu)=\mu\lambda\ln\lambda - (1+\mu\lambda)\ln\left(\frac{1+\mu\lambda}{1+\mu}\right)$
for $\lambda,\mu>0$ admits the following approximations:
\begin{enumerate}
\item For all $\mu > 0$ and $\lambda \in (0, 1]$, $H(\lambda, \mu) \geq  \frac{\mu}{2(1 + \mu)} \cdot (\lambda - 1)^2$;
\item For all $\mu > 0$ and $\lambda \in [1, 2]$, $H(\lambda, \mu) \geq  \frac{\mu}{4(1 + \mu)} \cdot (\lambda - 1)^2$;
\item For all $\mu > 0$ and $\lambda \geq 2$, $H(\lambda, \mu) \geq \frac{\mu}{4(1 + \mu)} \cdot (\lambda - 1)$;
\item For all $0 < \mu \leq \frac{1}{3}$ and $\lambda \geq 3$, $H(\lambda, \mu) \geq \frac{\mu\lambda}{4} \cdot \ln \left(\min\left\{\lambda, \frac{1}{\mu}\right\}\right)$;
\item For all $0 < \lambda \leq \frac{1}{3}$ and $\mu \geq 3$, $H(\lambda, \mu) \geq \frac{1}{4} \cdot \ln\left( \min\left\{\frac{1}{\lambda}, \mu\right\}\right)$.
\end{enumerate}
\label{prop:approx}
\end{proposition}



\subsection{Comparison with \cite{janson2018tail}}

In \cite{janson2018tail} a slightly different geometric random variable $X'=X+1$ was considered;
hence the tail bounds have to be carefully converted under the context of our geometric random variable model.

More specifically, Theorems~2.1 and 3.1 of \cite{janson2018tail} imply that
\begin{align*}
&\Pr\left[\bar X_n\geq\lambda\mu\right] \leq \exp\left\{-n\cdot ((\lambda - 1)(1-p) - \ln ((\lambda - 1)(1-p)+1))\right\}, & (\lambda>1);\\
&\Pr\left[\bar X_n\leq\lambda\mu\right] \leq \exp\left\{-n\cdot  ((\lambda - 1)(1-p) - \ln ((\lambda - 1)(1-p)+1))\right\}, & (\lambda<1);\\
\end{align*}
where $p$ is the geometric random variable parameter defined in Eq.~(\ref{eq:geometric}).
Because $z - \ln(z + 1) \leq z^2/2$ for all $z>0$, we have
\begin{equation}
 (\lambda - 1)(1-p) - \ln ((\lambda - 1)(1-p)+1) \leq \frac{\mu^2}{2(1+\mu)^2}\cdot (\lambda-1)^2.
\label{eq:janson}
\end{equation}

Comparing Eq.~(\ref{eq:janson}) with the simplified forms of $H(\lambda,\mu)$ in Proposition \ref{prop:approx},
we observe that Eq.~(\ref{eq:janson}) has an extra $\mu/(1+\mu)$ factor.
This means that when $\mu\to 0^+$ is very small, our results are tighter than \cite{janson2018tail}.

\subsection{Comparison with \cite{agrawal2017mnl}}

In Appendix D of \cite{agrawal2017mnl}, the following were established for sums of i.i.d.~geometrically distributed random variables:
for $\delta>0$, 
\begin{align*}
&\Pr\left[\bar X_n>(1+\delta)\mu\right] \leq \left\{\begin{array}{ll}\exp\left(-\frac{n\mu\delta^2}{2(1+\delta)(1+\mu)^2}\right)& \text{if }\mu\leq 1,\\
\exp\left(-\frac{n\mu^2\delta^2}{6(1+\mu)^2}\left(3-\frac{2\delta\mu}{1+\mu}\right)\right)& \text{if }\mu\geq 1, \delta\in(0,1);\\\end{array}
\right.\\
&\Pr\left[\bar X_n<(1-\delta)\mu\right] \leq \left\{\begin{array}{ll}\exp\left(-\frac{n\mu\delta^2}{6(1+\mu)^2}\left(3-\frac{2\delta\mu}{1+\mu}\right)\right)& \text{if }\mu\leq 1, \delta\in(0,1),\\
\exp\left(-\frac{n\mu^2\delta^2}{2(1+\mu)^2}\right)& \text{if }\mu\geq 1, \delta\in(0,1).\\\end{array}\right.
\end{align*}

We make several remarks comparing the above results from \cite{agrawal2017mnl} with our tail bounds.
First, the case of $\mu\geq 1$ and $\delta>1$ is missing in the above tail bounds, while our Theorem \ref{thm:main-positive} and Proposition \ref{prop:approx}
cover all cases of $\mu,\delta>0$.
In addition, in the case of $\mu\geq 3$ and $0<\lambda<1/3$, we have $\delta\geq 2/3$ and therefore 
$
\frac{n\mu^2\delta^2}{2(1+\mu)^2} \leq \frac{n}{2}.
$
On the other hand, by Proposition \ref{prop:approx} we know that $H(\lambda,\mu)\geq \frac{1}{4}\ln(\min\{\lambda^{-1},\mu\})$.
Hence our tail bound is sharper when $\mu\geq e^2$ and $\lambda<1/e^2$.
Finally, our forms of $H(\lambda,\mu)$ and its simplifications in Proposition \ref{prop:approx} are more user-friendly compared to \cite{agrawal2017mnl}.

\subsection{Proof of Proposition \ref{prop:approx}}

To simplify notations, we use $H_\lambda'(\cdot,\mu)$ and $H_\lambda''(\cdot,\mu)$ to denote $\partial H(t,\mu)/\partial t$
and $\partial H^2(t,\mu)/\partial t^2$, respectively.
Also define $H(0,\mu) := \lim_{\lambda\to 0^+}H(\lambda,\mu)$.

\noindent\textbf{\underline{Case 1: $\mu\geq 0$ and $\lambda\in(0,1]$.}}
Note that $H(1,\mu) = H_\lambda'(1,\mu)= 0$. 
By Taylor expansion with Lagrangian remainders, for any $\lambda\in(0,1]$ there exists $\lambda^*\in[\lambda,1]$ (depending on $\lambda$) such that
$
H(\lambda, \mu) = H_\lambda''(\lambda^*,\mu) (\lambda - 1)^2/2.
$
The first inequality of Proposition \ref{prop:approx} is then proved by noting that
$
H_\lambda''(\lambda^*,\mu) = \frac{\mu}{\lambda^*  (1 + \mu  \lambda^*)} \geq \frac{\mu}{1 + \mu}.
$

\noindent\textbf{\underline{Case 2: $\mu>0$ and $\lambda\in[1,2]$.}}
Let $\xi(\lambda) \eqdef H(\lambda, \mu) -  \frac{\mu}{4(1 + \mu)} \cdot (\lambda - 1)^2$. Because $\xi(1) = 0$, it suffices to show
$
\xi'(\lambda) = \mu  \ln\lambda - \mu  \ln\left(\frac{1 + \mu  \lambda}{1 + \mu}\right) -  \frac{\mu}{2(1 + \mu)} \cdot (\lambda - 1) \geq 0
$
for all $\lambda\in[1,2]$.
We know $\xi'(\lambda)$ is a concave function, because $\xi'''(\lambda) = \frac{\partial^3 H}{\partial\lambda^3}(\lambda, \mu) = - \frac{\mu  (1  + 2\mu  \lambda)}{\lambda^2  (1 +  \mu \lambda)^2} \leq 0$. Note also that $\xi'(1) = 0$. To settle our claim, it remains to show $\xi'(2) \geq 0$. Indeed,
$
\xi'(2) = -\mu  \ln\left(1 -  \frac{1}{2(1 + \mu)}\right) - \frac{1}{2} \cdot \frac{\mu}{1 + \mu} {\geq} 0,
$
where the last inequality holds because $-\ln(1 - z) \geq z$ for all $z \in (0, 1)$. This completes the proof of the second inequality in Proposition \ref{prop:approx}.

\noindent\textbf{\underline{Case 3: $\mu>0$ and $\lambda\geq 2$.}}
Because $-\ln(1 - z) \geq z$ for all $z \in (0, 1)$, 
\begin{equation}
\label{eq:4}
H_\lambda'(\lambda, \mu)
= -\mu  \ln\left(1 - \frac{1 - \lambda^{-1}}{1 + \mu}\right)
\geq \frac{1}{2} \cdot \frac{\mu}{1 + \mu}\;\;\;\;\forall \lambda\geq 2.
\end{equation}

By Taylor expansion with Lagrangian multiplier, for any $\lambda \geq 2$ there exists $\lambda^* \in [2, \lambda]$ depending on $\lambda$ such that
$
H(\lambda, \mu) = H(2, \mu) +H_\lambda'(\lambda*,\mu) \cdot (\lambda - 2)
\geq H(2, \mu) + \frac{\mu}{2(1+\mu)} \cdot (\lambda - 2),
$
where the last inequality follows by Eq.~(\ref{eq:4}).
Note also that $H(2, \mu) \geq  \frac{\mu}{4(1 + \mu)}$ by the second property of Proposition \ref{prop:approx}.
We then conclude that $H(\lambda, \mu) \geq \frac{1}{4} \cdot \frac{\mu}{1 + \mu} \cdot (2\lambda - 3)
\geq \frac{1}{4} \cdot \frac{\mu}{1 + \mu} \cdot (\lambda - 1)$, for all $\lambda \geq 2$.

\noindent\textbf{\underline{Case 4: $0<\mu\leq1/3$ and $\lambda\geq 3$.}}
It can be checked that for all $\lambda\geq 3$, 
\begin{equation}
\label{eq:5}
H\left(\lambda, {\lambda}^{-1}\right) = \ln\left[  \left(2 + \lambda + {\lambda}^{-1}\right)/4\right] \geq (1/4)\cdot \ln\lambda.
\end{equation}

When $\lambda \geq1/\mu \geq 3$, we also have $H_\lambda'(\lambda,\mu)
= \mu \cdot \ln\left(\frac{1 + \mu}{\lambda^{-1} + \mu}\right)
\geq \mu \cdot \ln\left(\frac{1 + \mu}{2\mu}\right)
\geq \frac{ \mu }{4} \ln\left(\frac{1}{\mu}\right)$.
Subsequently,
$
H(\lambda, \mu) -\frac{\mu\lambda}{4} \cdot \ln \left(\frac{1}{\mu}\right)
\geq H\left(\frac{1}{\mu}, \mu\right) - \frac{1}{4}  \ln \left(\frac{1}{\mu}\right)
{\geq} 0,
$
where the last inequality holds by Eq.~(\ref{eq:5}).

In the case of $1/\mu\geq\lambda\geq 3$, we know $\big[H(\lambda, \mu) - \frac{\mu\lambda}{4} \ln\lambda\big]$ is a concave function of $\mu$, because
$
\frac{\partial^2}{\partial \mu^2} \left[H(\lambda, \mu) - \frac{1}{4} \cdot \mu \lambda \cdot \ln\lambda\right] = \frac{\partial^2 H}{\partial \mu^2}(\lambda, \mu) = -\frac{(\lambda - 1)^2}{(1 + \mu)^2  (1 + \mu \lambda)} \leq 0.
$
By Eq.~(\ref{eq:5}) and the fact that $\lim_{t\to 0^+}H(\lambda, t) - \frac{1}{4}  t \lambda  \ln\lambda = 0$, we have
$H(\lambda, \mu) - \frac{1}{4} \cdot \mu \lambda \cdot \ln\lambda \geq 0$.
%
Combining both cases of $1/\mu\leq\lambda$ and $1/\mu\geq\lambda$ 
we complete the proof of the fourth property in Proposition \ref{prop:approx}.

\noindent\textbf{\underline{Case 5: $\mu\geq 3$ and $0<\lambda\leq1/3$.}}
It can be checked that for all $0<\lambda\leq 1/3$,
\begin{equation}
\label{eq:6}
H\left(\lambda, {\lambda}^{-1}\right) = \ln\left[ \left(2 + \lambda + {\lambda}^{-1}\right)/4\right] \geq 1/4\cdot \ln(\lambda^{-1}).
\end{equation}

Because $H_\lambda'(\lambda, \mu)
= -\mu \cdot \ln\left(1 + \frac{\lambda^{-1} - 1}{1 + \mu}\right) \leq 0$ for all $\lambda \in (0, 1]$, we have
forall $0 < \lambda \leq {1}/{\mu} \leq {1}/{3}$ that
$
H(\lambda, \mu) \geq H\left(\mu^{-1}, \mu\right){\geq} 1/4\cdot \ln\mu = 1/4 \cdot \ln\left( \min\left\{\lambda^{-1}, \mu\right\}\right),
$
where the second inequality holds by Eq.~(\ref{eq:6}).

On the other hand, for the case of $0 < {1}/{\mu} \leq \lambda \leq {1}/{3}$,
 because $\ln(1 + z) \leq z$ for all $z \in [0, +\infty)$, we have
 $
\frac{\partial H}{\partial \mu}(\lambda, \mu) = \frac{1 - \lambda}{1 + \mu} - \lambda \cdot \ln\left(1 + \frac{\lambda^{-1} - 1}{1 + \mu}\right) \geq 0
$
for all $\lambda\in(0,1]$.
Consequently,
$
H(\lambda, \mu) \geq H\left(\lambda, \lambda^{-1}\right){\geq} 1/4\cdot \ln\left(\lambda^{-1}\right) = 1/4 \cdot \ln\left(\min\left\{\lambda^{-1}, \mu\right\}\right),
$
where the second inequality follows Eq.~(\ref{eq:6}).
Combining both cases we complete the proof of the fifth property of Proposition \ref{prop:approx}.


\section{Proof of Theorem \ref{thm:main-positive}}\label{sec:proof-positive}


\noindent\textbf{\underline{Geometric distribution, upper tail ($\lambda>1$).}}
Using Chernoff bound, 
\[
\Pr\big[\bar X_n \geq \lambda \mu\big]
= \Pr\left[e^{t\bar X_n} \geq e^{t \lambda \mu}\right]
\leq \frac{1}{e^{t \lambda n \mu}} \cdot \prod\limits_{i = 1}^n \bE\left[e^{t X_i}\right], \;\;\;\; 0<t < -\ln(1 - p),
\]
where
$
\bE\left[e^{t X_i}\right] = {p}/[{1 - (1 - p) \cdot e^t}] = {1}/[{(1 + \mu) - \mu \cdot e^t}]
$
for $i=1,\cdots,n$.
Define $f(t) \eqdef \lambda\mu \cdot t + \ln\big[(1 + \mu) - \mu \cdot e^t\big]$. 
Note also that $-\ln(1-p)=\ln(1+\mu^{-1})$.
We then have
\begin{equation}
\label{eq:1}
\Pr\big[\bar X_n \geq \lambda \mu\big] \leq \exp\left\{-n \cdot f(t)\right\}, \;\;\;\;0<t<\ln(1+\mu^{-1}).
\end{equation}

For the function $f(\cdot)$ and its derivative $f'(\cdot)$, the following properties hold:
\begin{itemize}
\item $f'(t) = \lambda\mu - \frac{\mu \cdot e^t}{(1 + \mu) - \mu \cdot e^t}$ is a decreasing function when $0 < t < \ln\left(1 + \mu^{-1}\right)$;
\item $f'(t) = 0$ iff $\frac{\mu \cdot e^t}{(1 + \mu) - \mu \cdot e^t} = \lambda\mu$, or iff $t$ equals $T_1 \eqdef \ln\left(\frac{1 + \mu^{-1}}{1 + \mu^{-1} \cdot \lambda^{-1}}\right)$.
\end{itemize}
Because $0 < T_1 < \ln\big(1 + \mu^{-1}\big)$, $f(\cdot)$ attains its maximum at $f(T_1) = H(\lambda, \mu)$. 
Eq.~(\ref{eq:1}) then implies the desired upper tail bound.

\noindent\textbf{\underline{Geometric distribution, lower tail ($\lambda < 1$).}}
Again using the Chernoff bound,
\[
\Pr\big[\bar X_n \leq \lambda  \mu\big]
= \Pr\left[e^{-t\bar X_n} \geq e^{-t \lambda n \mu}\right]
\leq e^{t \lambda n \mu} \cdot \prod\limits_{i = 1}^n \bE\left[e^{-t X_i}\right], \;\;\;\; t>0,
\]
where $\bE\left[e^{-t X_i}\right]= {p}/[{1 - (1 - p) \cdot e^{-t}}] = {1}/[{(1 + \mu) - \mu \cdot e^{-t}}]$ for $i=1,\cdots,n$.
Define $g(t) \eqdef \lambda\mu \cdot t - \ln\big[(1 + \mu) - \mu \cdot e^{-t}\big]$. Then
\begin{equation}
\label{eq:2}
\Pr\big[\bar X_n \geq \lambda  \mu\big] \leq \exp\left\{n \cdot g(t)\right\}, \;\;\;\;t>0.
\end{equation}

For the function $g(\cdot)$ and its derivative $g'(\cdot)$, the following properties hold:
\begin{itemize}
\item $g'(t) = \lambda\mu - \frac{\mu \cdot e^{-t}}{(1 + \mu) - \mu \cdot e^{-t}}$ is an increasing function on $t \in (0, +\infty)$;
\item $g'(t) = 0$ iff $\frac{\mu \cdot e^{-t}}{(1 + \mu) - \mu \cdot e^{-t}} = \lambda\mu$, or iff $t$ equals $T_2 \eqdef \ln\left(\frac{1 + \mu^{-1} \cdot \lambda^{-1}}{1 + \mu^{-1}}\right)$.
\end{itemize}
Because $T_2 > 0$ for all $\lambda \in (0, 1)$, the minimum of $g$ is attained at $g(T_2) = -H(\lambda, \mu)$, which implies the desired lower tail bound.

\noindent\textbf{\underline{Exponential distribution, upper tail ($\lambda>1$).}}
By Chernoff bound, 
\[
\Pr\big[\bar Y_n \geq \lambda  \mu\big]
= \Pr\left[e^{t\bar Y_n} \geq e^{t \lambda  \mu}\right]
\leq \frac{1}{e^{t \lambda n \mu}} \cdot \prod\limits_{i = 1}^n \bE\left[e^{t Y_i}\right], \;\;\;\; 0<t<\mu^{-1},
\]
where $\bE\left[e^{t Y_i}\right] ={\int}_0^{+\infty}\mu^{-1} \cdot e^{-z/(\mu+tz)}  dz = {1}/[{1 - \mu  t}]$.
Define $h(t) \eqdef \lambda\mu t + \ln(1 - \mu  t)$.
We then have
\begin{equation}
\label{eq:3}
\Pr\big[\bar Y_n \geq \lambda  \mu\big] \leq \exp\left\{-n \cdot h(t)\right\}, \;\;\;\; 0<t<\mu^{-1}.
\end{equation}

For the function $h(\cdot)$ and its derivative $h'(\cdot)$, the following properties hold:
\begin{itemize}
\item $h'(t) = \lambda\mu - \frac{\mu}{1 - \mu \cdot t}$ is a decreasing function when $0 < t < \mu^{-1}$;
\item $h'(t) = 0$ iff $\frac{\mu}{1 - \mu \cdot t} = \lambda\mu$, or iff $t$ equals $T_3 \eqdef \mu^{-1} \cdot \left(1 - \lambda^{-1}\right)$.
\end{itemize}
Because $0 < T_3 < \mu^{-1}$, the maximum of $h(\cdot)$ is attained at $h(T_3)=G(\lambda)$,
which implies the desired upper tail bound.

\noindent\textbf{\underline{Exponential distribution, lower tail ($\lambda<1$).}}
This part is directly implied by Theorem~5.1(iii) of \cite{janson2018tail}.

\section{Proof of Corollary \ref{cor:negative}}\label{sec:proof-negative}

We shall prove the following theorems, 
which imply Corollary \ref{cor:negative}.
\begin{theorem}
Let $X_1,\cdots,X_n$ be geometrically distributed random variables; define $\Delta_U(\lambda,\mu):=7/6 + \ln(\lambda+2/\mu)$
and $\Delta_L(\lambda,\mu) := 1/6 + 3/2\cdot \ln(1+1/\lambda\mu)$.
For all $\lambda>1$ and all $n\geq 1$, 
\begin{equation}
\Pr\left[\bar X_n\geq \lambda\mu\right] \geq \exp\left\{-n\cdot\left[H(\lambda,\mu) + \frac{\ln(2\pi n)}{2n} + \frac{\Delta_U(\lambda,\mu)}{n}\right]\right\}.
\label{eq:negative-geometric-upper}
\end{equation}
In addition, for $\lambda \in(0,1)$ and all $n\geq 1/\lambda\mu$, 
\begin{equation}
\Pr\left[\bar X_n\leq \lambda\mu\right] \geq \exp\left\{-n\cdot\left[H(\lambda,\mu)+\frac{\ln(2\pi n)}{2n} + \frac{\Delta_L(\lambda,\mu)}{n}\right]\right\}.
\label{eq:negative-geometric-lower}
\end{equation}
\label{thm:negative-geometric}
\end{theorem}

\begin{theorem}
Let $Y_1,\cdots,Y_n$ be exponentially distributed random variables.
For all $\lambda\geq 1$ and $n\geq 1$, 
\begin{equation}
\Pr\left[\bar Y_n\geq\lambda\mu\right]\geq \exp\left\{-n\cdot\left[G(\lambda) + \frac{\ln(2\pi n)}{2n} + \frac{1}{12n^2}\right]\right\}.
\label{eq:negative-exponential-upper}
\end{equation}
In addition, for $\lambda\in(0,1)$ and $n\geq 1$,
\begin{equation}
\Pr\left[\bar Y_n\geq\lambda\mu\right] \geq \exp\left\{-n\cdot\left[G(\lambda)+\frac{\ln(2\pi n)}{2n} + \frac{1}{12n^2}\right]\right\}.
\label{eq:negative-exponential-lower}
\end{equation}
\label{thm:negative-exponential}
\end{theorem}

In the rest of this section we prove Theorems \ref{thm:negative-geometric} and \ref{thm:negative-exponential}.

\noindent\textbf{\underline{Geometric distribution, upper tail ($\lambda>1$).}}
We say a random variable $Z$ follows a binomial distribution with parameters $m\in\mathbb N$ and $p\in[0,1]$ if $\Pr[Z=k]=\binom{m}{k}p^k(1-p)^{m-k}$ for $k\in\{0,1,\cdots,m\}$.
Abbreviate by $B(m,1/(1+\mu))$ the binomial distribution parameterized by $m$ and $p=1/(1+\mu)$.
Given $\lambda \in [1, +\infty)$, the event $\{\bar X_n\geq \lambda \mu\}$ corresponds to fewer than $n$ successes in $\lceil\lambda n \mu\rceil + n - 1$ trials.
Formally,
\begin{align*}
\Pr\big[\bar X_n \geq \lambda  \mu\big]
= & \Pr\big[n\bar X_n \geq \lceil\lambda n \mu\rceil\big]
= \Pr_{Z\sim B(\lceil\lambda n\mu\rceil + n -1, 1/(1+\mu))}\left[Z < n\right]
& \\
= & \sum\limits_{k = 0}^{n - 1} \binom{\lceil\lambda n \mu\rceil + n - 1}{k}  \left(\frac{1}{1 + \mu}\right)^k  \left(\frac{\mu}{1 + \mu}\right)^{\lceil\lambda n \mu\rceil + n - 1 - k} \\
\geq & \binom{\lceil\lambda n \mu\rceil + n - 1}{n - 1}  \left(\frac{1}{1 + \mu}\right)^{n - 1}  \left(\frac{\mu}{1 + \mu}\right)^{\lceil\lambda n \mu\rceil} \\
\geq & \binom{\lceil\lambda n \mu\rceil + n}{n} \cdot \frac{n}{\lambda n \mu + 1 + n}  \left(\frac{1}{1 + \mu}\right)^{n - 1}  \left(\frac{\mu}{1 + \mu}\right)^{\lambda n \mu + 1}.
\end{align*}

Invoking Stirling's approximation (see \cite{robbins1955remark}, also summarized as Lemma \ref{lem:stirling} in the appendix), we have
\begin{flalign*}
\frac{1}{n}  \ln\binom{\lceil\lambda n \mu\rceil + n}{n}
\geq & \lambda \mu  \ln\left(1 + \frac{1}{\lambda \mu + 1 / n}\right) + \ln(1 + \lambda \mu) - \frac{\ln(2\pi n)}{2n} - \frac{1}{n} \cdot\frac{1}{6} \\
= & H_1(\lambda, \mu) - \frac{\ln(2\pi n)}{2n} - \lambda\mu  \ln\left(1 + \frac{1}{\lambda\mu}  \frac{1 / n}{1 + \lambda\mu + 1 / n}\right) - \frac{1}{6n} \\
\overset{(\dagger)}{\geq} & H_1(\lambda, \mu) - \frac{\ln(2\pi n)}{2n} - \frac{1/n}{1 + \lambda\mu + 1 / n} - \frac{1}{6n} \\
\geq & H_1(\lambda, \mu) - \frac{\ln(2\pi n)}{2n} - \frac{7}{6} \cdot \frac{1}{n},
\end{flalign*}
where $H_1(\lambda, \mu) \eqdef \mu  \lambda\ln\left(1 + {1}/({\mu  \lambda})\right) + \ln(1 + \mu  \lambda)$, 
$H_2(\lambda, \mu) \eqdef \ln(1 + \mu) + \lambda \mu\ln\left(1 + {1}/{\mu}\right)$,
and $(\dagger)$ holds because $\ln(1 + z) \leq z,\forall z \geq 0$.
Subsequently, 
\begin{flalign*}
\frac{n}{\lambda n \mu + 1 + n}
 \left(\frac{1}{1 + \mu}\right)^{n - 1}  \left(\frac{\mu}{1 + \mu}\right)^{\lambda n \mu + 1}
\geq & \exp\left\{-n \cdot \left[H_2(\lambda, \mu) + \frac{1}{n}  \ln\left(\lambda + \frac{2}{\mu}\right)\right]\right\}.
\end{flalign*}

This completes the proof of Eq.~(\ref{eq:negative-geometric-upper}) since  $H(\lambda, \mu) = H_2(\lambda, \mu) - H_1(\lambda, \mu)$.


\noindent\textbf{\underline{Geometric distribution, lower tail ($\lambda<1$).}}
Recall the definitions of $H_1(\lambda,\mu)=\mu\lambda\ln(1+1/(\mu\lambda))+\ln(1+\mu\lambda)$ and $H_2(\lambda,\mu)=\ln(1+\mu)+\lambda\mu\ln(1+1/\mu)$.
For $\lambda\in(0,1)$, we have
\begin{align*}
\Pr\big[\bar X_n \leq \lambda \mu\big]
= & \Pr\big[n\bar X_n \leq \lfloor\lambda n \mu\rfloor\big]
= \Pr_{Z\sim B(\lfloor\lambda n\mu\rfloor + n, 1/(1+\mu))}\left[Z \geq n\right]
& \\
= & \sum\limits_{k = n}^{\lfloor\lambda n \mu\rfloor + n} \binom{\lfloor\lambda n \mu\rfloor + n}{k}  \left(\frac{1}{1 + \mu}\right)^k  \left(\frac{\mu}{1 + \mu}\right)^{\lfloor\lambda n \mu\rfloor + n - k} \\
\geq & \binom{\lfloor\lambda n \mu\rfloor + n}{n}  \left(\frac{1}{1 + \mu}\right)^n  \left(\frac{\mu}{1 + \mu}\right)^{\lfloor\lambda n \mu\rfloor}.
\end{align*}

Note also that
$$
\left(\frac{1}{1 + \mu}\right)^n  \left(\frac{\mu}{1 + \mu}\right)^{\lfloor\lambda n \mu\rfloor} \geq \left(\frac{1}{1 + \mu}\right)^n  \left(\frac{\mu}{1 + \mu}\right)^{\lambda n \mu} = \exp\left\{-n \cdot H_2(\lambda, \mu)\right\}.
$$ 

With the condition $n \geq {1}/({\lambda\mu})$ which is equivalent to $\lfloor\lambda n \mu\rfloor \geq 1$, it follows from Lemma \ref{lem:stirling} (Stirling's approximation) that
\begin{flalign*}
\frac{1}{n}  \ln\binom{\lfloor\lambda n \mu\rfloor + n}{n}
\geq & \left(\lambda \mu - \frac{1}{2n}\right)  \ln\left(1 + \frac{1}{\lambda \mu}\right) + \ln\left(1 + \lambda \mu - \frac{1}{n}\right) - \frac{\ln(2\pi n)}{2n} - \frac{1}{6n} \\
= & H_1(\lambda, \mu) - \frac{\ln(2\pi n)}{2n} - \frac{1}{2n}  \ln\left(1 + \frac{1}{\lambda\mu}\right) + \ln\left(1 - \frac{1}{n} \cdot \frac{1}{1 + \lambda\mu}\right) - \frac{1}{6n} \\
\overset{(\ddagger)}{\geq} & H_1(\lambda, \mu) - \frac{\ln(2\pi n)}{2n} - \left[\frac{3}{2}  \ln\left(1 + \frac{1}{\lambda\mu}\right) + \frac{1}{6}\right]  \frac{1}{n},
\end{flalign*}
where $(\ddagger)$ holds because 
$
\left(1 - \frac{1}{n} \cdot \frac{1}{1 + \lambda\mu}\right)^n \geq 1 - \frac{1}{1 + \lambda\mu} = \left(1 + \frac{1}{\lambda\mu}\right)^{-1}
$.
This completes the proof of Eq.~(\ref{eq:negative-geometric-lower}) by noting that $H(\lambda, \mu) \equiv H_2(\lambda, \mu) - H_1(\lambda, \mu)$.

\noindent\textbf{\underline{Exponential distribution, upper tail ($\lambda>1$).}}
Consider a Poisson counting process $\{N(t):t\geq 0\}$ with rate $1/\mu$;
that is, $\Pr[N(t)=n]=(\lambda t)^ne^{-\lambda t}/n!$ for $n\in\mathbb N$. 
We then have for $\lambda>1$ that
\begin{flalign*}
\Pr\big[\bar Y_n \geq \lambda \mu\big]
= & \Pr\big[\bar Y_n > \lambda  \mu\big]
= \Pr\big[N(\lambda n \mu) < n\big]
& \text{($\bar Y_n$ has continuous density)} \\
= & \sum\limits_{k = 0}^{n - 1} \frac{(\lambda n)^k}{k!} \cdot e^{-\lambda n}
\geq \frac{(\lambda n)^{n - 1}}{(n - 1)!} \cdot e^{-\lambda  (n - 1)} \\
= & \frac{e^{\lambda}}{\lambda} \cdot \frac{(\lambda n)^n}{n!} \cdot e^{-\lambda n}
\geq \frac{(\lambda n)^n}{n!}  e^{-\lambda n}.
& \text{($e^z \geq z$ for all $z \geq 1$)}
\end{flalign*}

Note also that $\frac{(\lambda n)^n}{n!} \cdot e^{-\lambda n}$ can be re-written as $\exp\left\{-n \cdot \left[G(\lambda) + \frac{\ln n!}{n} - \ln n + 1\right]\right\}$. 
Using Lemma \ref{lem:stirling} (Stirling's approximation) we complete the proof of Eq.~(\ref{eq:negative-exponential-upper}).

\noindent\textbf{\underline{Exponential distribution, lower tail ($\lambda<1$).}}
Consider again the Poisson counting process $\{N(t):t\geq 0\}$ with rate ${1}/{\mu}$. For $\lambda \in (0, 1)$, we have
\begin{align*}
\Pr\big[\bar Y_n \leq \lambda \mu\big]
= & \Pr\big[N(\lambda n \mu) \geq n\big]
= \sum\limits_{k = n}^{\infty} \frac{(\lambda n)^k}{k!} \cdot e^{-\lambda n}
\geq \frac{(\lambda n)^n}{n!} \cdot e^{-\lambda n}. &
\end{align*}

Eq.~(\ref{eq:negative-exponential-lower}) can then be proved by following the same lines as in the proof of Eq.~(\ref{eq:negative-exponential-upper}) above.


\noindent\textbf{Stirling's approximation.}
The following lemma is from \cite{robbins1955remark}.
\begin{lemma}
For all $m\in\mathbb N$, $m>0$, it holds that
$$
\left(m + \frac{1}{2}\right)  \ln m - m + \frac{\ln(2\pi)}{2} + \frac{1}{12m + 1} \leq \ln (m!) \leq \left(m + \frac{1}{2}\right)  \ln m - m + \frac{\ln(2\pi)}{2} + \frac{1}{12m}.
$$
\label{lem:stirling}
\end{lemma}

\vspace{-0.2in}


\bibliographystyle{alpha}
\bibliography{refs}

\end{document}